# Closed Constant Curvature Space Curves
## Hermann Karcher, Bonn

Abstract: We use ODEs and symmetry arguments to construct closed constant curvature space curves, first on cylinders, next on tori, at last with the Frenet-Serret equations.

While I was still teaching I only knew closed constant curvature space curves which were pieced together from circle segments and helix segments. But smooth examples are easily accessible.

### 1. Examples on Cylinders

First roll the plane isometrically onto a cylinder of radius $R$:

$$F : \begin{pmatrix} x \\ y \end{pmatrix} \mapsto \begin{pmatrix} x \\ R\cos(y/R) \\ R\sin(y/R) \end{pmatrix}.$$

In the plane we describe a curve by its rotation angle against the x-axis, $\alpha(s) = \int_0^s \kappa_g(\sigma) d\sigma$, where $\kappa_g$ is the curvature of the plane curve, or its geodesic curvature when rolled onto the cylinder:

$$c'(s) := \begin{pmatrix} \cos(\alpha(s)) \\ \sin(\alpha(s)) \end{pmatrix}, \quad c(s) := \int_0^s c'(\sigma) d\sigma.$$

The cylinder has normal curvature 0 in the x-direction and $1/R$ in the y-direction. The space curvature $\kappa$ of $F \circ c$ is therefore given by

$$\kappa^2 = \sin^4(\alpha(s))/R^2 + \kappa_g^2(s) = \sin^4(\alpha(s))/R^2 + (\alpha'(s))^2.$$

This is a first order ODE for $\alpha(s)$, if we want $\kappa = const$.

This first order ODE is Lipschitz, if we look for curves with $\kappa > 1/R$:

$$\alpha'(s) = +\sqrt{\kappa^2 - \sin^4(\alpha(s))/R^2} > 0.$$

The solution curves are, in the plane, convex curves. They reach $\alpha = \pi/2$ in finite time. They are closed because the normals at $\alpha = 0$ and at $\alpha = \pi/2$ are lines of reflectional symmetry (of the plane curve).

For $\kappa \leq 1/R$ the ODE has some resemblence to the ODE $f' = \sqrt{1-f^2}$ of the sine function: it is not a Lipschitz ODE, with non-uniqueness along the constant solution $\alpha(s) = \arcsin(\sqrt{\kappa \cdot R})$. As with the sine-ODE we can differentiate the square of the ODE, cancel $2\alpha'(s)$ and obtain a second order Lipschitz ODE:
$$\alpha''(s) = -2\sin^3(\alpha(s))\cos(\alpha(s))/R^2.$$

If we choose $\kappa < 1/R$, then the second order ODE forces $\alpha'(s)$ to change sign when $\alpha(s)$ reaches $\alpha_{\max}$ given by $\sin^2(\alpha_{\max}) = \kappa \cdot R < 1$. The solution curves oscillate around a parallel to the x-axis and look a bit like sin-curves.

If we choose $\kappa = 1/R$, then $\alpha_{\max} = \pi/2$. We see that the circles $\alpha(s) := \pi/2$ are solutions of the second order ODE. By uniqueness, no solution which starts with $\alpha(0) < \pi/2$ can reach $\pi/2$ in finite time, it has to converge to $\pi/2$ asymptotically. The corresponding curve $F \circ c$ therefore spirals towards one of the circle-latitudes of the cylinder!

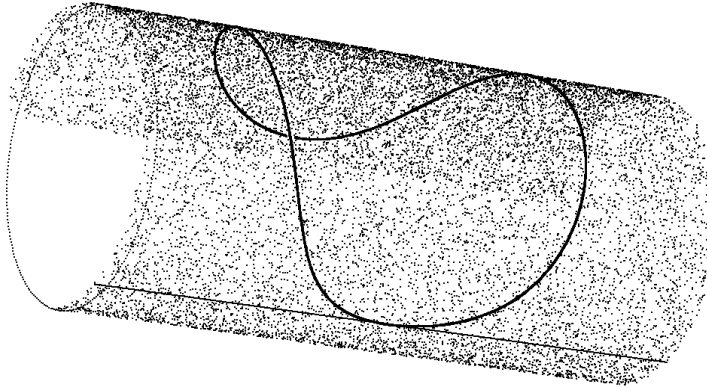

Convex curve in the plane, rolled onto a cylinder to a constant curvature space curve.

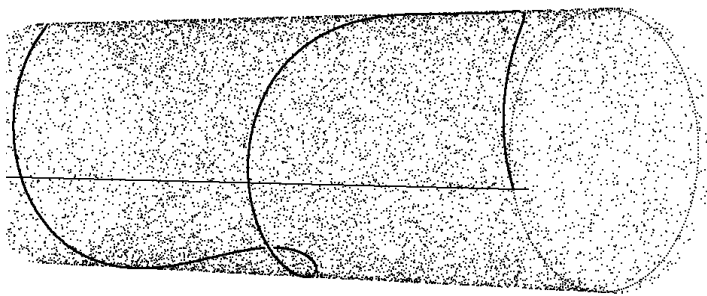

Periodic curve on cylinder, a constant curvature space curve. Geodesic curvature changes sign when crossing the drawn line.

## 2. Examples on Tori

The idea is the same as on the cylinder and works for many surfaces of revolution which also have a symmetry plane orthogonal to the rotation axis. A curve on the surface has constant space curvature $\kappa$ if its geodesic curvature $\kappa_g$ and its normal curvature $\kappa_n$ satisfy $\kappa^2 = \kappa_g^2 + \kappa_n^2$. Since the normal curvature depends only on the tangent of the curve there is again the easy case where we choose $\kappa^2 > \max \kappa_n^2$ and compute $\kappa_g(c'(s)) = +\sqrt{\kappa^2 - \kappa_n(c')^2}$ to get a second order ODE for the curve:

$$c''(s) = \kappa_n(c'(s)) \cdot N(c(s)) + \kappa_g(c'(s)) \cdot (c'(s) \times N(c(s))).$$

We start the integration on the equator, direction vertically up. Because the geodesic curvature is bounded away from zero, the curve will turn until it meets some meridian orthogonally. Reflection in the plane of this meridian and reflection in the equator plane complete the initial quarter arc to a closed curve. On the torus one can start at the inner or the outer equator.

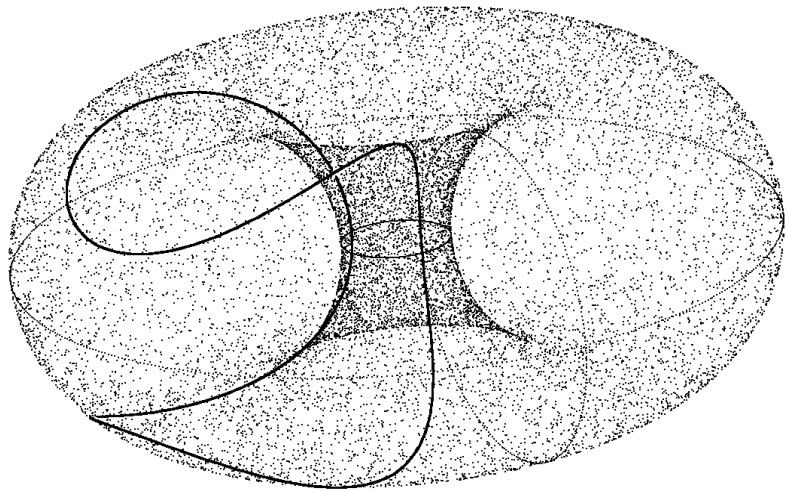

Constant curvature space curve, symmetric to the equator plane. The geodesic curvature is $\kappa_g \geq 0$.

With more effort we can find constant curvature space curves which oscillate around the equator. The geodesic curvature therefore needs to change sign, when the curve crosses the equator. That means, we need to choose as the space curvature the normal curvature of the torus in the direction $c'$ at the point where it crosses the equator. Therefore the integration of the ODE starts on the equator, the initial direction $c'(0)$ is a free parameter, the space curvature is computed as $\kappa = \mathrm{abs}(\kappa_n(c'(0)))$.

Up from the equator the normal curvature decrease and therefore the geodesic curvature increases. Again the angle of the curve against the meridians increases until it meets some meridian orthogonally. Reflection in the meridian plane continues the curve back to the equator and 180 degree rotation around the torus normal gives the next half-wave. In general these curves do not close. One needs to adjust eiither $c'(0)$ or the size of the torus until an even number of the initial half waves fits just once around. The curve closes smoothly, because the final tangent has the same angle with the equator as the initial tangent.

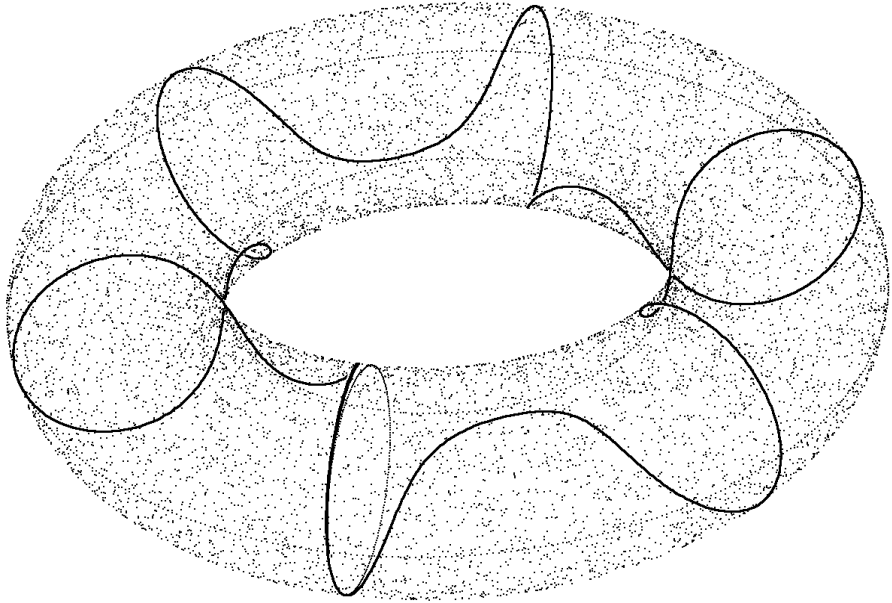

Constant curvature space curve, oscillating around the outer equator. The geodesic curvature changes sign where the curve crosses the equator.

### 3. Examples via the Frenet Equations

A closed space curve of constant curvature $\kappa(s) = \kappa$ must have a periodic torsion function, for example

$$\tau(s) = c \cdot \sin(s) + d \cdot \sin(2s) + e \cdot \sin(3s).$$

*The Frenet-Serret Equations:*

$$\begin{aligned} e_1' &= \kappa(s) \cdot e_2(s) \\ e_2' &= -\kappa(s) \cdot e_1(s) - \tau(s) \cdot e_3(s) \\ e_3' &= \tau(s) \cdot e_2(s) \\ c'(s) &= e_1(s) \end{aligned}$$

then determine the space curve $c$ uniquely from $c(0), e_1(0), e_2(0), e_3(0)$. We employ symmetries to find closed curves.

### 3.1 Reflection Symmetry in Normal Planes

A normal plane at $c(s_1)$ is a symmetry plane, iff the torsion function is skew symmetric w.r.t $s_1: \tau(s_1 - s) = -\tau(s_1 + s)$.

With our choice of $\tau$ the symmetry planes occur at $s = n \cdot \pi$, $n \in \mathbb{Z}$. Reflections in two neighboring symmetry planes generate all other symmetry planes so that the symmetry planes all intersect in one line. If the angle between two neighboring symmetry planes is a rational multiple of $\pi$, then the curve will be closed. In other words: in almost any 1-parameter family of curves there will be closed ones. Many of them will have selfintersections. The embedded ones look somewhat similar to those examples above which oscillate around the torus equator.

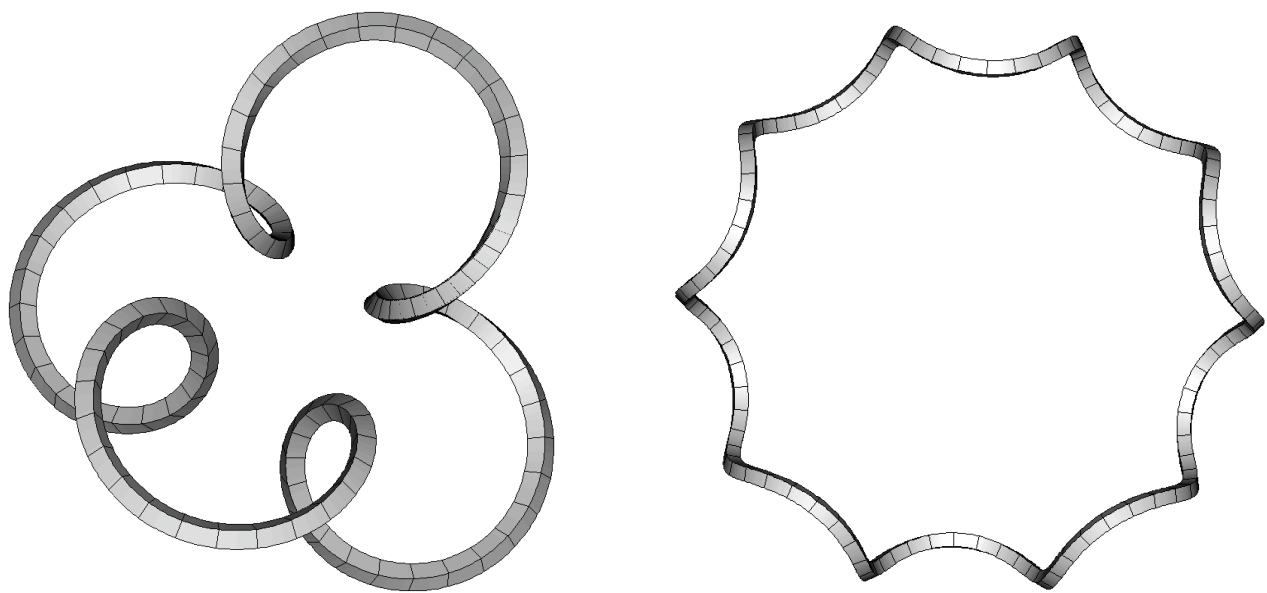

Constant curvature space curves with $\tau(s) = c \cdot \sin(s)$. Angle between neighboring normal symmetry planes: $90°$(left) and $36°$(right). Anaglyph 3D-images can be viewed at:
`virtualmathmuseum.org/SpaceCurves/constant_curvature/constant_curvature.html`.

### 3.2 $180°$-Rotation Symmetry around Principal Curvature Normals

We choose
$$\tau(s) = b + c \cdot \sin(s) + e \cdot \sin(3s).$$

Then $\tau(s)$ is symmetric w.r.t. $s_1 = \pi/2 + n \cdot \pi$, $n \in \mathbb{Z}$, that is: $\tau(s_1 - s) = \tau(s_1 + s)$. This implies that the curve $c$ is symmetric w.r.t. $180°$ rotation around the principal curvature normal $e_2(s_1)$. Clearly, if two such neighboring symmetry normals intersect, then all of them intersect in the same point. And if in addition the angle between two neighbors is a rational multiple of $\pi$, then the solution $c$ of the Frenet-Serret equations is a closed space curve of constant curvature. To find such curves, a 2-parameter problem has to be solved.

Next an observation comes in which reduces the problem to applications of the intermediate value theorem: The distance between neighboring symmetry normals depends in a surprisingly regular way on the constant coefficient $b$ of the Fourier polynomial $\tau(s)$. This can be used - for any choice of the other parameters $c, e$ - to adjust $b$ such that the symmetry normals all pass through one point. Here "surprisingly regular" means that this adjustment of $b$ causes no numerical difficulties. With this adjustment of $b$ always assumed, we consider the angle between neighboring symmetry lines as a function of the other parameters $c, e$. Again, in non constant 1-parameter families of solution curves there will be closed ones. Moreover, since one only rarely hits a local extremum of the angle function, one can find solution curves with nicely chosen angles such as $\pi/2, \pi/3, 2\pi/3, \pi/4, \pi/5...$ etc.

Since these curves are not constructed on well known surfaces, it is not so easy to see their shape in space from a printed image. We refer again to the anaglyph 3D-images at:

virtualmathmuseum.org/SpaceCurves/constant_curvature/constant_curvature.html.

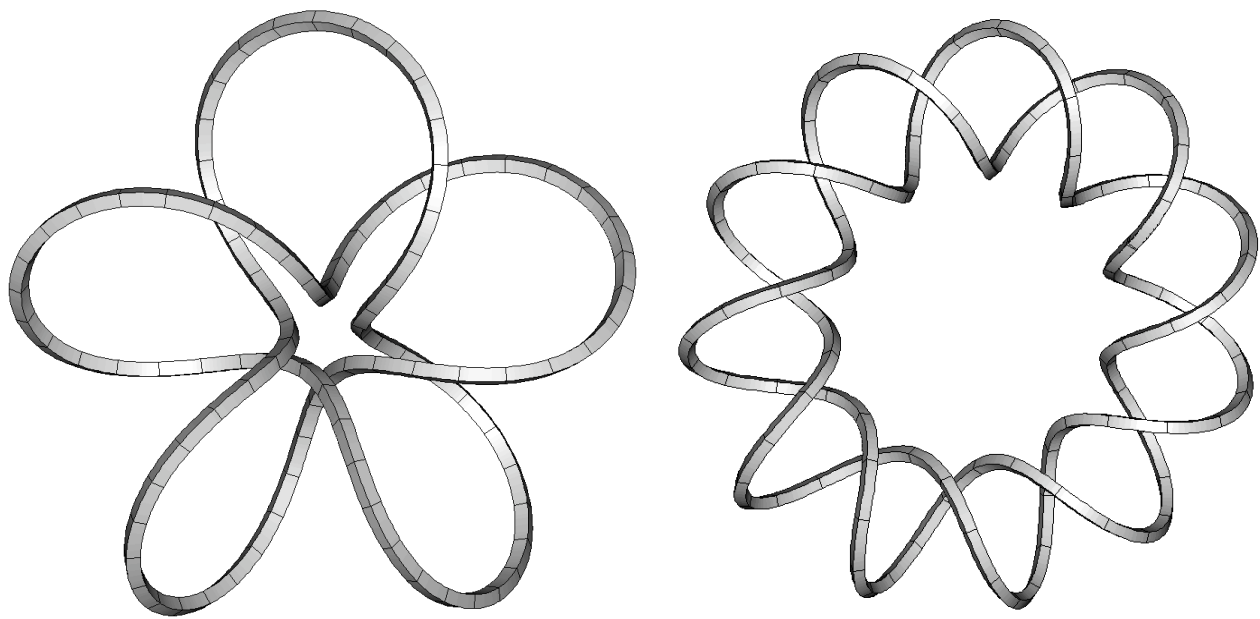

A 5-2-knot and a 11-2-knot of constant curvature, with $\tau(s) > 0$ for the second one. The principal normals at the outer most points are the symmetry normals. They meet the curve again at the inner most points.

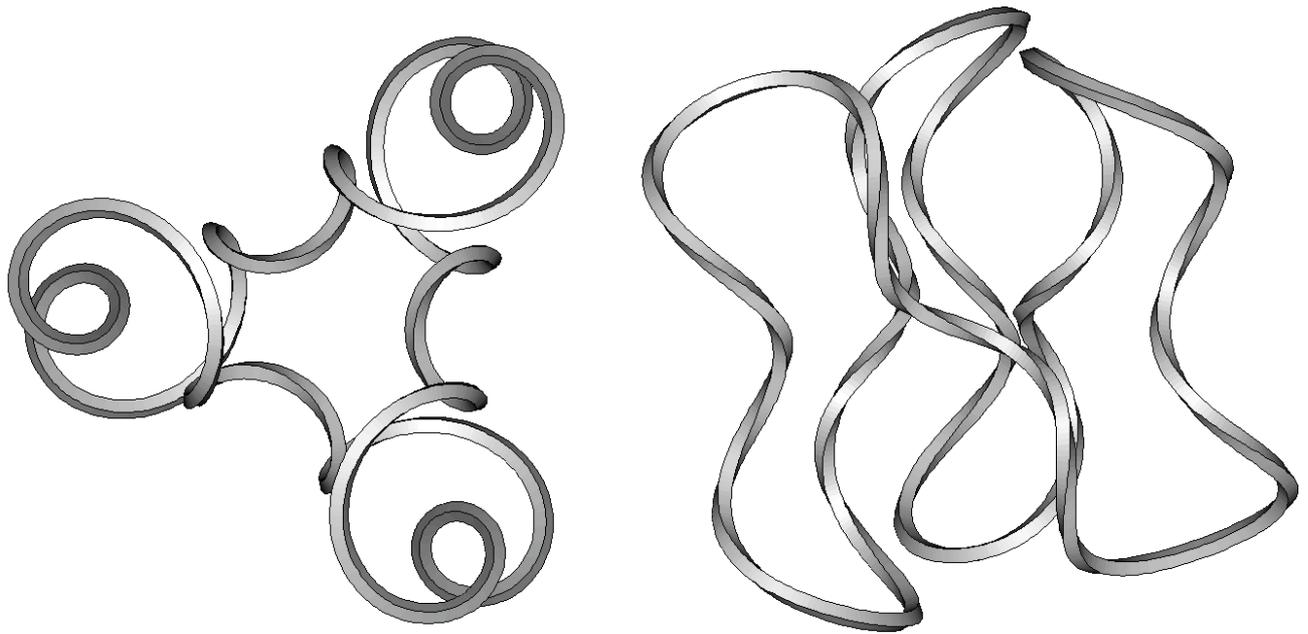

Top and side view of a more complicated space curve of constant curvature.

## Bibliography

I am not aware of publications which use ODEs from differential geometry to construct closed space curves of constant curvature.

Hermann Karcher
retired from
Dept. Math. Bonn University
unm416@uni-bonn.de